\documentclass[12pt]{amsart}
\usepackage{amsmath,amsthm,amsfonts,amssymb}
\begin{document} 
\newcommand{\tensor}{\otimes}
\newcommand{\B}{{\mathbb B}}
\newcommand{\C}{{\mathbb C}}
\newcommand{\N}{{\mathbb N}}
\newcommand{\Q}{{\mathbb Q}}
\newcommand{\Z}{{\mathbb Z}}
\renewcommand{\P}{{\mathbb P}}
\newcommand{\R}{{\mathbb R}}
\newcommand{\rc}{\subset}
\newcommand{\rank}{\mathop{rank}}
\newcommand{\trace}{\mathop{tr}}
\newcommand{\dimc}{\mathop{dim}_{\C}}
\newcommand{\Lie}{\mathop{Lie}}
\newcommand{\Ad}{\mathop{Ad}}
\newcommand{\Auto}{\mathop{{\rm Aut}_{\mathcal O}}}
\newcommand{\alg}[1]{{\mathbf #1}}
\newtheorem{definition}{Definition}
\newtheorem*{claim}{Claim}
\newtheorem{corollary}{Corollary}
\newtheorem*{Conjecture}{Conjecture}
\newtheorem*{SpecAss}{Special Assumptions}
\newtheorem{example}{Example}
\newtheorem*{remark}{Remark}
\newtheorem*{observation}{Observation}
\newtheorem{fact}{Fact}
\newtheorem*{remarks}{Remarks}
\newtheorem{lemma}{Lemma}
\newtheorem{proposition}{Proposition}
\newtheorem*{theorem}{Theorem}
\title{%
Dense Random Finitely Generated Subgroups of Lie Groups
}
\author {J\"org Winkelmann}
\address{%
J\"org Winkelmann \\
 Institut Elie Cartan (Math\'ematiques)\\
 Universit\'e Henri Poincar\'e Nancy 1\\
 B.P. 239\\
 F-54506 Vand\oe uvre-les-Nancy Cedex\\
 France
}
\email{jwinkel@member.ams.org\newline\indent{\itshape Webpage: }%
http://www.math.unibas.ch/\~{ }winkel/
}
\thanks{
{\em Acknowledgement.}
The author wants to thank the organizers of ``Groups and
Probability'', L.~Pyber and M.~Ab\'ert, for the invitation
to participate in this conference and E.~Breuillard for
asking the question which is answered by this article.
}%
\maketitle
\section{The result}
The purpose of this paper is to prove the following theorem:
\begin{theorem}
Let $G$ be a connected real Lie group of dimension $n$.

Then there exists a relatively compact open neighbourhood $W$ of $e$
in $G$ such that for $n+1$  randomly chosen elements
$g_0,\ldots,g_n$ the generated subgroup will be dense in $G$
with probability one.
\end{theorem}

\section{Preparations}
\subsection{The probability measure}
Every locally compact group admits a ``Haar measure'' which is
finite iff the group is compact.%
\footnote{Here it does not matter whether we consider
left- or rightinvariant Haar measures.}
 The measure of a compact subset
is always finite, the measure of an open subset always non-zero.
Therefore, by renormalizing we obtain an induced probability
measure for every relatively compact open subset.

In other words, if, given a Lie group $G$ and a relatively compact
open subset $W\subset G$, we say ``$k+1$ randomly chosen elements 
$g_0,\ldots,g_k$ in $W$
have the property $X$ with probability one'', then this means the following:

{\em Let 
\[
\Sigma
=\{(g_0,\ldots,g_k)\in W^{k+1}: 
\text{$(g_0,\ldots,g_k)$ has property $X$}\}.
\]
Let $\mu$ be a Haar measure on $G$ and $\mu^{k+1}$ the induced
product measure on $G^{k+1}$. Then $\mu^{k+1}(\Sigma)=\mu^{k+1}(W^{k+1})$.}

\subsection{Zassenhaus neighbourhoods}
Inspired by the classical notion of a {\em Zassenhaus neighborhood}
we introduce:
\begin{definition}
Let $G$ be a connected real Lie group. For an element $g\in G$
let $\zeta_g:G\to G$ be the commutator map, i.e.,
$\zeta_g(h)=ghg^{-1}h^{-1}$.

A $Z$-neighbourhood
$W$ is an open relatively compact subset of $G$ such that
$e\in W$ and $\lim_{k\to\infty}(\zeta_g)^k(x)=e$
for all $g,x\in W$.
\end{definition}
Developing the group law on $G$ into a power series, it becomes
obvious that every real Lie group admits a $Z$-neighbourhood.

We use this definition rather than the classical notion of
a Zassenhaus neighborhood%
\footnote{Classically, a neighbourhood $U$ of $e$ in $G$ is
called a {\em Zassenhaus neighbourhood} iff every  discrete subgroup of $G$
generated by its intersection with $U$ must be nilpotent.}
  because for our definition it is clear that it is stable
under taking quotients of Lie groups:
\begin{fact}\label{fact-Z}
Let $\phi:G\to H$ be a surjective homomorphism of real
Lie groups and $W$ a $Z$-neighbourhood in $G$.

Then $\phi(W)$ is a $Z$-neighbourhood in $H$.
\end{fact}

\subsection{Regular elements}

We recall the notions of Cartan subalgebras and regular elements in
Lie groups.
As standard references we use 
\cite{Bou}, \cite{Dix} and \cite{Enz}.

\begin{definition}
Let $\alg g$ be a finite-dimensional Lie algebra over a field $k$.
A Lie subalgebra $\alg h$ is called {\em Cartan subalgebra}
if it is nilpotent and equals its own normalizer
(i.e. $[x,a]\in\alg h\,\,\forall a\in \alg h$ implies $x\in\alg h$.)
\end{definition}

\begin{definition}
Let $G$ be a Lie group. An element $g\in G$ is called
{\em regular}
if the multiplicity of $1$ as root of the characteristic
polynomial of $Ad(g)$ is minimal.
\end{definition}
This definition implies immediately that the set of all non regular
elements in a connected Lie group constitutes a nowhere dense real analytic
subset. In particular (see e.g.~\cite{jwtop} for more details):
\begin{fact}\label{almost-all-regular}
Let $G$ be a connected real Lie group and $G_{reg}$ the subset of
all regular elements.

Then $G\setminus G_{reg}$ has Haar measure zero.
\end{fact}

\begin{lemma}\label{reg-norm-nil}
Let $G$ be a connected real Lie group, $g$ a regular element and
$H$ a connected closed normal subgroup of $G$ which contains $g$.

Then $G/H$ is nilpotent.
\end{lemma}
\begin{proof}
For a complex number $\lambda$ let $V_\lambda$ be the associated
weight space for $Ad(g)$, i.e.
\[
V_\lambda=
\left\{ v \in (\Lie G)\tensor_\R\C: \left(Ad(g)-\lambda I\right)^N(v)=0
\ \exists N>0\right\}.
\]
Then evidently $V_\lambda\subset\left(\Lie H\right)\tensor_\R\C$
for all $\lambda\ne 0$.
Therefore $V_0$ surjects onto $\left ( \Lie G/\Lie H\right)\tensor_\R\C$.
But for a regular element $g$ the weight space $V_0$ is a Cartan subalgebra
of $\Lie G$ and therefore a nilpotent Lie subalgebra of $\Lie G$.
Hence $G/H$ must be nilpotent.
\end{proof}
\subsection{The commutative case}
\begin{lemma}\label{ab-qu}
Let $H$ be a subgroup of $G=\left(\R^n,+\right)$.

Then either $H$ is dense in $G$ or there exists a 
surjective $\R$-linear map
$F:\R^n\to\R$ with $F(H)\subset\Z$.
\end{lemma}
\begin{proof}
Let $\bar H$ be the closure of $H$ in $G$.
The quotient $G/\bar H$ is a commutative connected real Lie group
and therefore isomorphic to some $(S^1)^d\times\R^g$.
Both $\R$ and $S^1$ admit surjective Lie group homomorphisms onto
$S^1$. Thus either $H$ is dense in $G$ or there is a surjective
Lie group homomorphism from $G/\bar H\to S^1$. Using the isomorphism
$S^1\simeq\R/\Z$ it is clear that in the latter case there is
a surjective $\R$-linear map $F$ from $\R^n$ to $\R$ with $F(H)\subset\Z$.
\end{proof}
\begin{lemma}\label{abelian}
Let $G=\R^n$. Then there is a set $\Sigma$ of Haar (i.e.~Lebesgue)
measure zero in $G^{n+1}$ such that for every
$(g_1,\ldots,g_{n+1})\in G^{n+1}\setminus\Sigma$
the subgroup of $G$ generated by the $g_i$ is dense.
\end{lemma}
\begin{proof}
Let $\Sigma_0$ denote the set of all $(g_1,\ldots,g_{n+1})$ for which
$G$ is not generated as a $\R$-vector space by $g_1,\ldots,g_n$.
For each $(g_1,\ldots,g_{n+1})\not\in\Sigma_0$ the element $g_{n+1}$
can be expressed as $\R$-linear combination of the other $g_i$:
\begin{equation}\label{eq}
g_{n+1}=\sum_{j=1}^n a_jg_j\quad \exists a_j\in\R
\end{equation}
We define $\Sigma_1\subset G^n\setminus\Sigma_0$
as the set of all $(g_1,\ldots,g_{n+1})$
such that the real numbers $1,a_1,\ldots, a_{n}$
are not $\Q$-linearly independent.

We set $\Sigma=\Sigma_0\cup\Sigma_1$.
Let $(g_1,\ldots,g_{n+1})\in G^{n+1}\setminus\Sigma$.
Let $(a_j)_{1\le j\le n}$ as in eq.~\ref{eq}.
Then (due to the definition of $\Sigma$) the real
numbers $1,a_1,\ldots,a_n$ are $\Q$-linearly independent.
As a consequence we obtain:

\[
\sum_j a_jb_j\not\in\Q
\]

for any choice of $b_1,\ldots, b_n\in\R$ unless all the $b_j$
are zero.

This implies: If $F:\R^n\to\R$ is a non-zero $\R$-linear map
with $F(g_1),\ldots,F(g_n)\in\Q$, then $F(g_{n+1})\not\in\Q$.

Therefore there does not exist a non-zero $\R$-linear map $F:G\to\R$
with $F(\Gamma)\subset\Q$ and a forteriori no such map
with $F(\Gamma)\subset\Z$.
By lemma~\ref{ab-qu} it thus follows that $\Gamma$ is dense
in $G$ whenever $(g_1,\ldots,g_{n+1})\not\in\Sigma$.
\end{proof}

\subsection{The nilpotent case}
We recall the notion of the (descending) central series:

Given a group $G$ and subsets $A,B\subset G$ let $[A,B]$ denote
the set of all commutators $aba^{-1}b^{-1}$ with $a\in A$ and
$b\in B$.

The groups $G^{[k]}$ of the central series are defined
recursively:
$G^1=G$ and if $G^{[k]}$ is already defined, one defines
$G^{[k+1]}$ as the closure of the subgroup of $G$ generated by the
elements of $[G,G^{[k]}]$.

A group $G$ is nilpotent iff there exists a number $n\in\N$
such that $G^{[n]}=\{e\}$.

By definition $G^{[2]}$ coincides with the
{\em commutator group $G'$}.

\begin{lemma}\label{nil-reduce}
Let $G$ be a real nilpotent Lie group
and $H$ a subgroup.
Assume that $HG^{[2]}$ is dense in $G$.

Then $H$ is dense in $G$.
\end{lemma}
\begin{proof}
We start the proof of
the lemma by claiming that $H^{[k]}G^{[k+1]}$ is dense in $G^{[k]}$ for all $k\in\N$.
We will show this claim
by induction on $k$.
Thus let us assume
that $H^{[k]}G^{[k+1]}$ is dense in $G^{[k]}$.
We have to show that under this assumption 
$H^{[k+1]}G^{[k+2]}$ must be dense in $G^{[k+1]}$.
Now  $G^{[k+1]}$ is topologically generated
by $[G,G^{[k]}]$. Since  $HG^{[2]}$ is dense in $G$ 
and $H^{[k]}G^{[k+1]}$ is dense in
$G^{[k]}$, it follows that the image of $[H,H^{[k]}]$ must be dense in
$[G,G^{[k]}]/A$ where $A$ 
is topologically
generated by $[G^{[2]},G^{[k]}]$ and $[G,G^{[k+1]}]$.
Since $[G^{[2]},G^{[k]}]$ and $[G,G^{[k+1]}]$ are both contained 
in $G^{[k+2]}$,
it follows that $H^{[k+1]}G^{[k+2]}$ is dense in $G^{[k+1]}$.
This proves the claim.

Next we show by induction in the opposite direction
that $H^{[k]}$ is dense $G^{[k]}$ for all $k\in\N$.
For $k$ sufficiently large this is trivially true since 
$H^{[k]}=G^{[k]}=\{e\}$ for all sufficiently large $k$.
Now let us assume that $H^{[k]}$ is dense in $G^{[k]}$.
Then $H^{[k-1]}=H^{[k-1]}\cdot H^{[k]}$ is dense
in $H^{[k-1]}\cdot G^{[k]}$ which in turn is dense
in $G^{[k-1]}$.
Thus we deduce by (descending) induction on $k$ that
$H^{[k]}$ is dense in $G^{[k]}$ for all $k\in \N$.
In particular, $H=H^{[1]}$ is dense in $G=G^{[1]}$.
\end{proof}

\begin{proposition}\label{prop-nil}
Let $G$ be a connected nilpotent real Lie group,
and $d=\dim(G/G')$.

Then there exists a subset $\Sigma\subset G^{d+1}$ of Haar measure zero
such that for every $(g_0,\ldots,g_d)\in G^{d+1}\setminus\Sigma$
the subgroup of $G$ generated by the $g_i$ is dense in $G$.
\end{proposition}

\begin{proof}
Let $\tau:G\to G/G'$ be the natural projection.
Thanks to lemma~\ref{abelian} we know that there is a subset
$\Sigma_0\subset (G/G')^{d+1}$ of Haar measure zero such that
$d+1$ elements $g_0,\ldots,g_d\in G/G'$ generate a dense subgroup
unless $(g_0,\ldots, g_d)\in\Sigma_0$.

Let $\Sigma=\tau^{-1}(\Sigma_0)$.
Then $\Sigma$ has Haar measure zero and for each
$(g_0,\ldots, g_d)\in G^{d+1}\setminus\Sigma$ the generated subgroup
$\Gamma$ has a dense image in $G/G'$.

Finally due to lemma \ref{nil-reduce} the condition
$\overline{\tau(\Gamma)}=G/G'$ implies $\overline{\Gamma}=G$.
\end{proof}

\section{Proof of the main result}

\begin{proposition}\label{prop-gen}
Let $G$ be a connected real Lie group of dimension $n$,
$W$ a $Z$-neighbourhood and
$0\le k \le n$.

Let $g_0,\ldots,g_k\in W$ be randomly chosen elements
and $H_k$ the closure of the subgroup of $G$ generated
by $g_0,\ldots,g_k$.

Then with probability one at least one of the following two
conditions is fulfilled:
\begin{itemize}
\item There is a normal closed Lie subgroup $M$ in $G$ with
$M\subset H_k$ such that $G/M$ is nilpotent.
\item $\dim H_k\ge k$.
\end{itemize}
\end{proposition}
\begin{proof}
The proof works by induction on $k$. For $k=0$ there is nothing to
prove.

Thus let us assume $k>0$ and let us furthermore assume that $H_{k-1}^0$
has the desired property. In other words, we assume that
\[
(g_0,\ldots,g_{k-1})\in W
\] 
are given such that
the closure $H_{k-1}$ of the 
group generated by these elements has the desired
property and we want to show that under
this assumption a randomly chosen element $g_k\in W$
has with probability one the property that $g_0,\ldots,g_{k-1}$
and $g_k$ together generate a subgroup with the desired property.

Since the set of non-regular elements of $G$ is a set of
Haar measure zero (fact~\ref{almost-all-regular}), 
we may and do assume that all the
$g_i$ are regular elements.

If $H_{k-1}^0$ contains a closed subgroup $M$ such that $H_{k-1}^0$
and $G/M$ nilpotent, so does $H_k^0$ and we are done.
Therefore we may assume that $H_{k-1}^0$ satisfies the second
property, i.e. $\dim H_{k-1}\ge k-1$.

We distinguish three cases:
\begin{enumerate}
\item $H_{k-1}^0$ is not normal in $G$.
\item $H_{k-1}^0$ is normal and $g_0\not\in H_{k-1}^0$.
\item $H_{k-1}^0$ is normal and $g_0\in H_{k-1}^0$.
\end{enumerate}

In the first case the normalizer of $H_{k-1}^0$ in $G$ is a proper
real Lie subgroup of $G$, implying that with probability one
$g_kH_{k-1}^0g_k^{-1}$ is not contained in $H_{k-1}^0$
for a randomly chosen $g_k\in W$.
But then $H_{k-1}^0$ and $g_kH_{k-1}^0g_k^{-1}$ generate a Lie subgroup
$I$ with $\dim(I)>\dim(H_{k-1})\ge k-1$ and $I\subset H_k$.
Hence $\dim H_k\ge k$ with probability one in this case.

In the second case let us consider the projection
$\pi:G\to A=G/H_{k-1}^0$.
Note that $\pi(W)$ is a $Z$-neighbourhood in $A$ (see fact~\ref{fact-Z}).
We are done if $A$ is nilpotent.
Thus we may and do assume that $A$ is not nilpotent.
Since we are in the second case, $\pi(g_0)\in\pi(W)\setminus\{e_A\}$.
Recall that we assumed $g_0$ to be a regular element in $G$.
Thus the Lie subalgebra of $\Lie G$ defined by the zero weight
space of $\Ad(g_0)$ is a Cartan subalgebra and in particular it is
nilpotent.
Since $A$ is not nilpotent,
 this Cartan subalgebra can not map 
surjectively on $A$. Therefore there is a non-zero complex
number $\lambda$ which occurs as eigen value for
$Ad(\pi(g_0))\in GL(\Lie A)$. As a consequence 
$\zeta_{\pi(g_0)}^N:A\to A$
is never a constant
map (here $\zeta_{\pi(g_0)}$ is the commutator map,
i.e.~$\zeta_{\pi(g)}(h)=ghg^{-1}h^{-1}$.)
Thus $\zeta_{\pi(g_0)}^N(\pi(g_k))\ne e$ with probability one for
all $N\in\N$ and a randomly chosen element $g_k\in W$.
On the other hand $\lim_{N\to\infty}\zeta^N_{\pi(g_0)}(\pi(g))=e$
for all $g\in W$,
because $\pi(W)$ is a $Z$-neighbourhood in $A$.
It follows that with probability
one the subgroup of $A$ generated by $\pi(g_0)$ and $\pi(g_k)$
is not discrete. This implies $\dim(H_k)>\dim(H_{k-1})$ and
thereby $\dim(H_k)\ge k$.

It remains to discuss the third case.
However, since $g_0$ is regular, the assumptions
``$H_{k-1}^0$ is normal in $G$ and $g_0\in H_{k-1}^0$''
do imply that $G/H_{k-1}^0$ is nilpotent (see lemma~\ref{reg-norm-nil}).
Thus we found a connected Lie subgroup $M$ with $G/M$ nilpotent and 
$M\subset H_k$, namely $M=H_{k-1}^0$.
\end{proof}

\begin{corollary}\label{cor-gen}
Let $G$ be a real connected Lie group of dimension $n$ and $W$ a
$Z$-neighbourhood.

For randomly chosen elements $g_0,\ldots, g_n\in W$ let $H$
be the closure of the subgroup
generated by the $g_i$.

Then with probability one there exists a normal closed Lie subgroup
$M$ of $G$ such that $G/M$ is nilpotent and $M\subset H$.
\end{corollary}
\begin{proof}
By prop.~\ref{prop-gen}
with $k=n$ we know that there is such a Lie subgroup $M$ unless
$\dim(H)\ge n$. But $\dim(H)\ge n$ implies $\dim(H)=\dim(G)$
and therefore $H=G$. In this case $M=G$ has the desired properties.
\end{proof}

Now we can prove the theorem:
\begin{proof}
Let $N$ denote the intersection of all $\ker\phi$ where $\phi$
runs through all Lie group homomorphisms from $G$ to any nilpotent
Lie group.
Then $N$ is a normal closed Lie subgroup of $G$ and $G/N$ is nilpotent.
Moreover, if $M$ is any normal closed Lie subgroup of $G$
for which $G/M$ is nilpotent, then $N\subset M$.
Let $\tau:G\to G/N$ denote the natural projection.

Now prop.~\ref{prop-nil} implies:

{\em If $W$ is a relatively compact open subset of $G$,
and $g_0,\ldots, g_n$ are randomly chosen elements in $W$,
then with probability one the elements $\tau(g_i)$ do generate
a dense subgroup of $G/N$.}

On the other hand, from cor.~\ref{cor-gen} we infer:

{\em If $W$ is a $Z$-neighbourhood in $G$, and
$g_0,\ldots,g_n$ are randomly chosen elements of $W$,
and $H$ is the closure of the subgroup of $G$ generated by the
$g_i$, then with probability one there is a normal closed
Lie subgroup $M$ of $G$ such that $G/M$ is nilpotent and
$M\subset H$. 
Since $G/M$ being nilpotent implies $N\subset M$,
we obtain that $N\subset H$ with probability one.%
}

Combined, these assertions prove the theorem:
If $\tau(H)$ is dense in $G/N$ and $N\subset H$, then
$H$ must be dense in $G$.
\end{proof}
\section{Optimality}
The theorem can not be improved on the number of generators
for an arbitrary Lie group, since a dense subgroup of $(\R^d,+)$
can not be generated by less than $d+1$ elements.

The relatively compact open subset $W$ in the statement
of the theorem can be any $Z$-neighborhood, but not an arbitrary
relatively compact open neighborhood:

\begin{proposition}
Let $n$ be a natural number and $G$ a connected real Lie group
which is not amenable.

Then there exists a relatively compact open subset $U$ in $G$
such that $n$ randomly chosen elements in $U$ generate a
non-dense subgroup with probability at least $\frac{n!}{n^n}$.
\end{proposition}
\begin{proof}
By the results of \cite{jwtop} we can find open subsets $V_1,\ldots,V_n$
in $G$ such that for any choice of $g_i\in V_i$ the subgroup of $G$
generated by the $g_i$ is discrete (hence not dense).
By shrinking the $V_i$ we may assume that they are all relatively
compact, are pairwise disjoint and
have the same volume with respect to the
Haar measure $\mu$ of $G$.
Define $U=\cup_iV_i$. Observe that $\mu(V_i)=\frac{1}{n}\mu(U)$
for each $i$ by our assumptions.

Let $\sigma$ be a permutation of $\{1,\ldots,n\}$.
If $g_i\in V_{\sigma(i)}$ for all
$i\in\{1,\ldots,n\}$, then the subgroup generated by the $g_i$
is discrete. For each $i$ the probability for $g_i\in V_{\sigma(i)}$ to hold
is $\frac{1}{n}$. The cardinality of all such permutations $\sigma$
is $n!$. Hence with probability at least $\frac{n!}{n^n}$ randomly chosen
elements $g_1,\ldots,g_n$ in $U$ do generate a discrete
(and therefore non-dense) subgroup of $G$.
\end{proof}
\section{Nilpotent Lie groups}
Above (see prop.~\ref{prop-nil}) we proved:

{\em Let $G$ be a connected nilpotent Lie group and $d=\dim(G/G')$.
Then there exists a subset $\Sigma\subset G^{d+1}$ of Haar measure
zero such that for any $(g_1,\ldots,g_{d+1})\in G^{d+1}\setminus
\Sigma$ the subgroup generated by the $g_i$ is dense in $G$.}

On the other hand, in \cite{jwtop} we proved:

{\em Let $G$ be a connected nilpotent Lie group and $d=\dim(G/G')$.
Then there exists a natural number $N\le d$ such that for all
$k\le N$ there exists a subset $\Sigma_k\subset G^k$ of Haar measure
zero such that for any $(g_1,\ldots,g_{k})\in G^k\setminus
\Sigma_k$ the subgroup generated by the $g_i$ is discrete in $G$.}

We would like to emphasize that it may happen that $N<d$. In
this case for all $k\in\{N+1,\ldots,d\}$ any $k$ randomly
chosen elements generate a subgroup which is neither dense
nor discrete (with probability one).

\begin{example}
Let $G$ be the simply-connected real nilpotent Lie group
associated to the Lie algebra $\left<A,B,C,D\right>$ with
$[A,B]=C$ and $[A,C]=D$ as the only non-trivial commutator relations
among the base vectors.

$G$ can also be described as follows: Take $\R^4$ as manifold and
define the group law via
\begin{multline*}
(a_1,b_1,c_1,d_1)\cdot (a_2,b_2,c_2,d_2)
=\\
\left(a_1+a_2,b_1+b_2,c_1+c_2+a_1b_2,
d_1+d_2+a_1c_2+\frac{a_1^2}{2}b_2 \right)
\end{multline*}

Let $Z$ be the center, i.e.
$Z=\{(0,0,0,d):d\in\R\}$
and $\Lie Z=\left< D \right>$.

Now let us chose two elements $g_i=(a_i,b_i,c_i,d_i)$
(for $i\in\{1,2\}$).
If $(a_1,b_1)$ and $(a_2,b_2)$ are $\R$-linearly independent elements
of $\R^2$ and $a_1,a_2$ are $\Q$-linearly independent elements of $\R$,
then the connected component of the closure of the 
subgroup $\Gamma$ of $G$ generated by
$g_1$ and $g_2$ equals $Z$. In particular, in this case
$\Gamma$ is neither discrete nor dense.
\end{example}

\section{Perfect Lie groups}
A Lie group $G$ is called {\em perfect} if it equals its own
commutator group. For these special class of Lie groups
better results follow immediately from the work
of Breuillard and Gelander. In particular, their results imply
the following:

{\em
Let $G$ be a perfect Lie group and let $k$ be a natural number such
that the Lie algebra of $G$ can be generated 
(as a Lie algebra) by $k$ elements.

Then there is a relatively compact
open neighbourhood $W$ of $e$ in $G$ such that
with probability one $k$ randomly chosen elements of $W$ will
generate a dense subgroup of $G$.
}

For a semisimple Lie group this implies that $2$ elements are
enough (because every semisimple Lie algebra is generated by
two elements), 
while for an arbitrary perfect Lie group it is at least
clear that $\dim G$ elements suffice.

Since our method needs $\dim G+1$ generators, the results
of Breuillard and Gelander yield a better result
for perfect Lie groups.

\end{document}